\documentclass[a4paper]{article}
\usepackage{amsmath}
\usepackage{amsthm}
\usepackage{amssymb}

\usepackage[latin1]{inputenc}

\usepackage[dvips=true,%
a4paper=true,%
pdftitle={On some bound and scattering states associated with the cosine kernel},%
pdfauthor={Jean-Francois Burnol},%
pdfstartview=FitH,%
pdfpagemode=None%
]{hyperref}

\def\ladate{January 3, 2008}

\newtheorem{theorem}{Theorem}

\newtheorem{proposition}[theorem]{Proposition}

\theoremstyle{definition}
\newtheorem{definition}[theorem]{Definition}

 \newcommand{\NN}{{\mathbb N}}

 \newcommand{\RR}{{\mathbb R}}
 \newcommand{\CC}{{\mathbb C}}
 \newcommand{\HH}{{\mathbb H}}

\newcommand{\cA}{{\mathcal A}}
\newcommand{\cB}{{\mathcal B}}
\newcommand{\cC}{{\mathcal C}}

\newcommand{\cE}{{\mathcal E}}
\newcommand{\cF}{{\mathcal F}}
\newcommand{\cG}{{\mathcal G}}

\newcommand{\cZ}{{\mathcal Z}}

\newcommand{\Un}{{\mathbf 1}}

\let\wh=\widehat
\let\wt=\widetilde

\begin{document}

\title{On some bound and scattering states associated with the cosine kernel}

\author{Jean-Fran\c cois Burnol}

\date{\ladate}

\maketitle

\begin{center}
\begin{small}
\noindent\begin{minipage}{4cm}
I.H.\'E.S.\\
Le Bois-Marie\\
35, route de Chartres\\
F-91440 Bures-sur-Yvette\\
France\\
burnol@ihes.fr\\
\end{minipage}\\

\emph{on leave from:}
Universit\'e Lille 1,
UFR de Math\'ematiques, 
F-59655 Villeneuve d'Ascq, 
France.\\
burnol@math.univ-lille1.fr

\begin{abstract}   It is explained how to provide self-adjoint
operators having scattering states forming a multiplicity one continuum and
bound states whose corresponding eigenvalues have an asymptotic density
equivalent to the one of the zeros of the Riemann zeta function. It is shown
how this can be put into an integro-differential form of a type recently
considered by Sierra.
\end{abstract}
\medskip

\end{small}
\end{center}

%%%%\pagebreak

%% \begin{small}

%% \tableofcontents

%% \end{small}

\setlength{\normalbaselineskip}{16pt}
\baselineskip\normalbaselineskip
\setlength{\parskip}{6pt}

\section{Introduction}

The present publication was motivated by a recent investigation of Sierra
\cite{Sierra4} (also \cite{Sierra2,Sierra3}). It will be explained how the
differential system and isometric expansion we obtained in
\cite{cras2003}, and the associated self-adjoint operators, provide, when moved
back to the ``$x$-picture''  a particular  example of a type of
differential-integral equation considered in  \cite[section IV]{Sierra4}. 

In \cite{cras2002} we determined ``explicitely'' certain objects depending on
a parameter $a>0$ and associated to the cosine (or sine) kernel. In
\cite{cras2003} we obtained the differential system (which turned out to be of
the Dirac type with coefficients involving some Fredholm determinants of the
cosine (or sine) kernel) obeyed  by the objects  from \cite{cras2002} as
functions of $\log(a)$. We mentioned the associated isometric spectral
expansion and also explained how this investigation had led to the realization
of the Fourier transform as a scattering, an objective which had resulted from
our definition and study of the ``$\log|x|+\log|p|$'' operator
(\cite{conductor}). 

In \cite{hankel} we gave a detailed exposition of these results. Some general
aspects  could be predicted in advance as they involve some Hilbert spaces of
entire functions in the sense of \cite{Bra}. The specific spaces had actually
already be defined in \cite{bra64,rov1}, among spaces  associated to general Hankel
transforms. The results from \cite{bra64,rov1} regarded mainly the spaces
associated to the Bessel function $J_0$. In \cite{cras2002} we obtained for
the cosine and sine kernels the ``explicit'' form taken by some functions
whose existence was a consequence of the general theory of \cite{Bra}: in
particular certain entire functions $\cA_a(s)$ and $\cB_a(s)$ having all their
zeros $\rho$ on a line, and such that the quotients $\frac{\cA_a(s)}{s-\rho}$
(resp. $\frac{\cB_a(s)}{s-\rho}$) give orthogonal bases of some Hilbert space
$K_a$. It was then shown in \cite{cras2003} that these objects obey
differential equations whose coefficients are expressed in terms of certain
Fredholm determinants. This observation also applied to the $J_0$ kernel and
we gave the detailed exposition of these results in \cite{hankel}, where the
chapters V and VI with minor adaptations apply to all Hankel kernels, in
particular to the cosine and sine kernels. A further advance had been realized
in the theory of the $J_0$ kernel and associated spaces whose nature was
revealed to be special relativistic \cite{einstein}. Other Hilbert
spaces of entire functions, directly associated with the Riemann zeta function
and which had not been considered before, were defined by the author in
\cite{cras2001} and shortly thereafter understood to be in connection with the
``co-Poisson'' formula.

In \cite{Sierra4}, Sierra considers integro-differential operators which are
singular rank two perturbations of $-i(x\frac d{dx}+\frac12)$ on
$L^2(0,+\infty;dx)$. We will explain how the objects, differential equations,
and spectral expansion from \cite{cras2002,cras2003,hankel}, when seen in the
``$x$ picture'' lead to the kind of equations considered by Sierra.  As was
mentioned in a number of our other papers (such as \cite{twosys}) the theory
leads to ``bound states'' associated to zeros on the critical line which have
an asymptotic density equivalent to the one applying to  the zeros of the
Riemann zeta function. These bound states are created on an interval
$(\log(a),+\infty)$ by a Dirichlet condition at $\log(a)$ and a potential
function which has an exponential increase at $+\infty$. The Dirac system seen
on $(-\infty,\log(a))$  leads to a purely continuous spectrum, because the
potential vanishes exponentially quickly at $-\infty$. Hence we do have here
the combination of scattering and bound states mentioned in section IV of
\cite{Sierra4} as a possible property of the integro-differential equations
considered there.  It is not clear if the particular example which is
considered in section VI of \cite{Sierra4}  truly leads to bound states, or
rather only to resonances.  On the other hand, the Dirac and Schrödinger
equations from \cite{cras2003} do have, as we mentioned in our publications,
this property of leading to ``quantum zeros'' lying on the critical line with
the expected density (a general density result has been given in
\cite{twosys}, see Theorem 7.7 and Remark 19).

%% We also include some results, interesting in their own sake and related to the
%% others mentioned here, although they are not explicitely involved in the
%% sections devoted to the comparison with \cite{Sierra4}. 

We shall not enter
here into other topics considered by Sierra, such as obtaining operators
leading exactly to the zeros of Riemann.

\clearpage

\section{Some notations}

For the sake of facilitating comparisons with \cite{Sierra4} we shall often
write $s=\frac12 + iE$ when using complex numbers although the letter $E$ is
also used with a completely different meaning. Also, our scalar products
$(\cdot\,|\,\cdot)$ will be conjugate linear in the left entry (bra) and
complex linear in the right entry (ket). As we refer to somewhat lengthy
developments from our previous papers we did not modify our notations too
much, in particular with respect to the use of the letters $A$ and $B$. We
will not here attempt to analyse our final equations in the spirit of
\cite[IV]{Sierra4} but only want to show how they emerge from our framework, so
the incompatibility of notations will remain virtual.

The Mellin transform is given by the
formula $\wh f(\frac12 + iE) = \int_0^\infty f(x)\frac{x^{-iE}}{\sqrt x}\,dx =
\int_0^\infty f(x) x^{-s}\,dx$ (and not with $+iE$) because this represents
the scalar product $(\psi_E^0|f)$, of $f$ with the
generalized eigenvector $\psi_E^0(x) = \frac{x^{iE}}{\sqrt x}$ of the ``free''
operator $H_0 = -i(x\frac d{dx} + \frac12)$, $H_0(\psi_E^0) = E\psi_E^0$. The
inverse Mellin transform is given by the formula $f(x) = \frac1{2\pi}
\int_{-\infty}^{+\infty} \wh f(\frac12 + iE)\frac{x^{iE}}{\sqrt x}\,dE =
\frac1{2\pi} \int_{\Re(s)=\frac12} \wh f(s) x^{s-1}\,|ds|$.

%% , which represents the expansion $\int (\psi_E^0|f)
%% \psi_E^0\,\frac{dE}{2\pi}$ of $f$ in the ``basis'' of the
%% $\psi_E^0$'s, $(\psi_{E}^0|\psi_{E'}^0) = 2\pi \delta(E-E')$.

%% It is sometimes useful to write $f(E)$ rather than the more explicit but long
%% $\wh f(\frac12 + iE)$, while maintaining the notation $f(x)$ for $f$ in the
%% $x$-picture. If we were to switch to the $q = \log(x)$ picture, and apply
%% similar abuse of notations we would then write $f(q) = \sqrt{x}f(x)$ because
%% this makes $\int_\RR |f(q)|^2\,dq = \int_0^\infty |f(x)|^2\,dx$ and $\frac
%% d{dq} = x\frac d{dx} + \frac12$. As arguably $f(q) = \sqrt{x}f(x)$ can be
%% confusing the use of the $q$-picture in fact more or less requires new letters
%% such that, for example,  $F(q) = \sqrt{x}f(x)$ or $\alpha(q)=
%% \sqrt{x}\psi_\alpha(x)$. In fact we will try to avoid the ``$q$-picture''
%% altogether.

%% When dealing with ``scattering states'' $\psi_E^{\text{scat}}$ which behave as
%% $C_0(E) x^{iE-\frac12}$ as $x\to0^+$ and as $C_\infty(E) x^{iE-\frac12}$ as
%% $x\to+\infty$ we will make a global choice of $C_0(E)$ which loses none of
%% those scattering states, for example $C_0(E) = 1$, or some other function of
%% $E$ which does not vanish for $E$ real.

%% We may sometimes say that $\wh f(s)$ is ``even'' or ``odd'', with the
%% meaning that $\pi^{-\frac s2}\Gamma(\frac s2) \wh f(s)$ is even or odd
%% under $s\mapsto 1-s$, which means that $f$ is invariant, respectively
%% skew-invariant, under the cosine transform $\cC$.

Regarding the Fourier transform, it will be used in the form of the cosine
transform   $\cC(f)(y) = \int_0^{+\infty} 2\cos(2\pi xy)f(x)dx$, for functions
or distributions supported in $[0,+\infty)$. The cosine kernel $2\cos(2\pi xy)$ has only the two eigenvalues $+1$
and $-1$ on $(0,\infty)$. When restricted to a finite interval $(0,a)$ it has
a discrete spectrum $1>|\lambda_1(a)|\geq |\lambda_2(a)| \geq \dots$ which
accumulates at zero. This restriction is conjugated to, hence has the same
spectrum as, the kernel $C_a(x,y) = 2a \cos(2\pi a^2 xy)$ acting on the
interval $(0,1)$. We will make use of the Fredholm determinants $\det(1 \pm
C_a)$. Some of their properties (in particular the asymptotic as
$a\to+\infty$) can be derived from knowledge of the corresponding facts for
$\det(1 -C_a^2)$. As the square $C_a^2$  is the Dirichlet kernel
$\frac{\sin(2\pi a^2 (x-y))}{\pi (x-y)}$ acting on the even subspace of
$L^2(-1,+1;dx)$, the question is reduced to the finite Dirichlet kernel, whose
properties  have been extensively studied in the literature in connection with
random matrices: see \cite{mehta} and the further references included in our
papers \cite{cras2003,hankel}.

\section{A support condition and associated Hilbert spaces}

To each square integrable function $f$ on $(0,+\infty)$, one can associate its
Mellin transform $\wh f(s) = \int_0^\infty f(x) x^{-s} dx$ on the critical
line $\Re(s) = \frac12$. This Mellin-Plancherel transform is a unitary
identification of $L^2(0,+\infty; dx)$ with $L^2(\Re(s)=\frac12;
\frac{|ds|}{2\pi})$. The closed subspace $L^2(1,+\infty; dx)$ (resp. $L^2(0,1;
dx)$) becomes identified with (boundary values of the functions of) the Hardy
space of the half-plane $\Re(s)>\frac12$ (resp. $\Re(s)<\frac12$). Let us
recall how the cosine transform $g(y) = \cC(f)(y) = \int_0^{+\infty}
2\cos(2\pi xy)f(x)dx$ is represented on the Mellin side. The formula is:
\begin{equation}
 \wh{\cC(f)}(s) = \chi(s) \wh f(1-s)\qquad\text{where}\qquad \chi(s) =
\pi^{s-\frac12}\frac{\Gamma(\frac{1-s}2)}{\Gamma(\frac s2)}
\end{equation} This is the same
$\chi(s)$,\footnote{let us observe in passing $\chi(s)\chi(1-s) = 1$ and using
formally $f(x)=\delta(x-1)$ that $\chi(s) = \int_0^\infty 2\cos(2\pi x)
x^{-s}dx$. This is valid as a semi-convergent integral, and also as an
identity of tempered distributions in $\log(x)$ and $\Im(s)$
($\Re(s)=\frac12$).} of modulus one on the critical line, which appears in the
functional equation $\zeta(s) = \chi(s)\zeta(1-s)$ of the Riemann zeta
function. Indeed this functional equation can be seen to be a manifestation of
the self-invariance under the cosine transform of the distribution
$\sum_{n\geq1} \delta(x-n) - \Un_{x>0}(x)$, whose Mellin transform is nothing
else but $\zeta(s)$ itself.

%%
%%  (alternatively the
%% self-invariance of the Dirac comb $\sum_{n\in\ZZ} \delta(x-n)$ under the 
%% Fourier transform $\int_\RR e^{i 2\pi xy} f(x) dx$ on the full line) 
%%

Of course $\zeta(s)$ is not square integrable on the critical line but the
functions $s\mapsto \frac{\zeta(s)}{s-\rho}$ are, with $\rho$  a zero of the
zeta function. Let us consider the vectors $\zeta_{\rho,k}(s) =
\frac{\zeta(s)}{(s-\rho)^k}$ where $\rho$ is a non-trivial zero of zeta, and
$k$ is an integer between $1$ and the multiplicity of $\rho$. Taking inverse
Mellin transforms we also tacitly consider the $\zeta_{\rho,k}$ to be vectors in
$L^2(0,+\infty;dx)$. The following holds  (\cite{twosys}):

\begin{theorem} The vectors $\zeta_{\rho,k}$  are linearly independent in the
sense that none is in the closure of the linear combinations of the
others. Furthermore the smallest closed subspace of $L^2(0,+\infty;dx)$
containing these vectors consists exactly of those functions $f$ in $L^2$
which are constant on $(0,1)$ and whose cosine transforms are also constant on
$(0,1)$.
\end{theorem}

We mentioned this theorem as motivation to study the space $K_1$ of functions
$f$ in $L^2$ which vanish identically on $(0,1)$ and whose cosine transform
also vanish identically on that same interval (the condition of vanishing is
slightly easier than the one of being constant). More generally we can
consider the space $K_a$ where the interval $(0,1)$ is replaced with the
interval $(0,a)$ (we would have been led
directly to these spaces if we had considered
the functions $\frac{L(s,\chi)}{s-\rho}$ where $\chi$ is an even Dirichlet
character):

\begin{definition} 
We let, for each $a>0$, $K_a$ be the subspace of $L^2(0,+\infty;dx)$
consisting of  functions which vanish
identically on $(0,a)$ and whose cosine transforms also vanish identically on
that same interval.
\end{definition}

The subspace $K_a$ is closed and its orthogonal complement is the
(non-orthogonal, but direct and closed) sum $L^2(0,a;dx) +
\cC(L^2(0,a;dx))$. 

\begin{proposition}
  If $f\in K_a$ then $\wh f(s)$ is an entire function which has trivial zeros at
  $s=0$, $-2$, $-4$, $-6$, \dots
\end{proposition}

Indeed $\wh f(s)$ as Mellin transform of $f$ which is supported in
$[a,+\infty)$ is analytic in the half-plane $\Re(s)>\frac12$, and its boundary
values on the line are the Mellin-Plancherel transform of $f$ (taken a priori
in the $L^2$-sense). The same applies to $g=\cC(f)$ and on the critical line
one has the almost everywhere identity $\wh g(s) = \chi(s)\wh{f}(1-s)$. So the
function $k(s) = \chi(s)\wh g(1-s)$ is analytic in the half-plane $\Re(s) <
\frac12$, and its boundary values from the left coincides with the boundary
values of $\wh f$ from the right. From general principles of complex analysis
(``edge-of-the-wedge'' theorems) the two analytic functions ``glue'' to give
an entire function. It would be easy to give in this one-dimensional case more
details for this glueing, but anyhow various simple elementary constructive
proofs of the Proposition are available, such as the one given in
\cite{cras2001}.\footnote{in \cite{cras2001} we  use $\int_0^\infty f(x)x^{s-1}\,dx$ whereas systematically in later references we always use
$\int_0^\infty f(x)x^{-s}\,dx$. Hence the $s$ from \cite{cras2001} is the $1-s$
from here.} The trivial zeros are a corollary of the identity $\wh{\cC(f)}(s)
= \chi(s)\wh f(1-s)$, because $\wh f(1-s)$ must have the zeros to compensate
the poles of $\chi(s)$ at $1$, $3$, $5$, \dots

The spaces $K_a$, mainly as spaces of the entire functions $\cF(s) =
\pi^{-\frac s2} \Gamma(\frac s2) \wh f(s)$ (but seen as functions of
the variable $z$ such that $s = \frac12 - iz$), had appeared earlier
in the literature: in de~Branges investigation \cite{bra64} where they
are studied as a special instance of spaces verifying the three
axiomatic properties\footnote{not needed here, see \cite{Bra}, and for
the $K_a$ spaces \cite{bra64,rov1} and \cite{cras2002,twosys,hankel}.}
whose general consequences are the subject of the book \cite{Bra}. The
well-known Paley-Wiener spaces (of Fourier transforms of $L^2$
functions of compact support) are another, simpler, instance of the
general theory of \cite{Bra}. A particularly interesting (and
challenging) aspect of the spaces $K_a$ is that the associated entire
functions are not of finite exponential type. The reader interested
generally speaking in \cite{Bra} is also referred to \cite{dymkean}
and to \cite{dym} and \cite{remling}, especially for the rather
intimate connection with standard aspects of spectral expansions
associated to differential equations of the Sturm-Liouville or
Schrödinger types.

In the process  in \cite{cras2002} and \cite{cras2003} of studying the
spaces $K_a$ directly  we have in particular recovered  by other means
the consequences predicted by the general theory of \cite{Bra}, in an
``explicit'' form. In this section and the next we briefly review some
of these results.

The evaluation maps $f\mapsto \wh f(s)$ are continuous linear forms on $K_a$,
they thus correspond to certain vectors in $K_a$  indexed by the complex
numbers. The space $K_a$ is ``elucidated'' once the mutual scalar products
between these evaluators are computed.  Let $\cZ_z$ be the entire function
such that $\cF(z) = (\cZ_z|\cF)$ for all completed Mellin transforms $\cF$ of
elements of $K_a$ (we may consider $\cZ_z$ either as an element of $K_a\subset
L^2(0,\infty;dx)$ or in its incarnation as an entire function $\cZ_z(s)$; the
dependency on $z$ is anti-analytic). Then, from \cite{Bra} it is known that
there must exist an entire function $\cE_a$ (not quite but almost unique) such
that:\footnote{we replace here the horizontal line always used in \cite{Bra}
by the critical line $\Re(s)=\frac12$. Hence formula \eqref{eq:repro1} does
not look as in \cite{Bra}.}
\begin{equation}\label{eq:repro1} \cZ_z(w) = (\cZ_w|\cZ_z) =
\frac{\overline{\cE_a(z)}\cE_a(w) - \cE_a(\overline
{1-z})\overline{\cE_a(\overline{1-w})}}{\overline z + w -1}
\end{equation}  
Such an $\cE_a$ once found also characterizes the space of the
entire functions $\cF(s)$ as follows: the ratio $\frac\cF{\cE_a}$ must belong
to the Hardy space of the half-plane $\Re(s)>\frac12$ and also the ratio
$\frac\cG{\cE_a}$ where $\cG(s) = \overline{\cF(\overline{1-s})}$. The formula
$\|f\|_2 = \|\frac\cF{\cE_a}\|_{\HH^2}$ holds, that is, somewhat surprisingly:
\begin{equation}
\frac1{2\pi} \int_\RR |\wh f(\frac12+iE)|^2\,dE = \frac1{2\pi}
\int_\RR
|{\cF}(\frac12+iE)|^2\,\frac{dE}{|\cE_a(\frac12+iE)|^2} =
\frac1{2\pi} \int_\RR \frac{|{\wh
f}(\frac12+iE)|^2}{|\wh{E_a}(\frac12+iE)|^2}\,{dE}
\end{equation}
We have introduced here the notation $\wh{E_a}(s) = \cE_a(s)/\pi^{-\frac
  s2}\Gamma(\frac s2)$ and indeed this proves useful when studying the spaces
$K_a$. 

The reader is referred to \cite{Bra} and  \cite{dymkean} for the general proof
(which is short) of existence of such ``$\cE$-functions'' when the
axioms of \cite{Bra} are verified. From the evaluator identity one deduces
that $|\cE_a(s)|^2 > |\cE_a(\overline{1-s})|^2$ when $\Re(s)>\frac12$, so all
zeros of $\cE_a(s)$ verify $\Re(s)\leq\frac12$. A zero on the symmetry line
would mean that all functions $\cF(s)$ vanish there. It can be shown that this
is not the case, so the zeros of $\cE_a(s)$ in fact verify
$\Re(s)<\frac12$. \footnote{Lagarias has considered in \cite{laghouches} an
entire function which has the properties of an $\cE$-function \`a la
de~Branges if and only if the Riemann hypothesis holds. Lagarias $\cE$
function would vanish at the locations on the critical line of the non-simple
zeros of zeta, if some exist.}

The author has contributed the ``explicit'' determination of a suitable
function $\cE_a$ \cite{cras2002}. It is realized as $\pi^{-\frac s2}
\Gamma(\frac s2) \wh{E_a}(s)$ where $E_a(x)$ is in fact a distribution
supported on $[a,+\infty)$ and whose cosine transform is also supported on
$[a,+\infty)$. To obtain $E_a$, we studied the vectors $X_s(x)$ which are such
that $\int_a^\infty f(x)X_s(x)\,dx = \wh f(s)$ for $f\in K_a$. For
$\Re(s)>\frac12$ the vector $X_s$ is simply the orthogonal projection to $K_a$
of $x^{-s}\Un_{x>a}(x)$. We determined directly the integrals $\int_a^\infty
X_z(x)X_w(x)\,dx$.  A completely detailed exposition with all proofs has been
given in our manuscript \cite{hankel}.
As explained there, the reader will
only need to replace  all occurences of the kernel $J_0(2\sqrt{xy})$ by the
kernel $2\cos(2\pi xy)$, and also to replace all occurences of $\Gamma(s)$ by
$\pi^{-\frac s2}\Gamma(\frac s2)$. The chapters of \cite{hankel} which are
related to the present publication are the fifth and the sixth. The
simplifications allowed by the kernel $J_0(2\sqrt{xy})$ are made use of only
in other chapters of \cite{hankel}.

%  
%% %% Let us hope that despite the referee's inspired 
%% %% comment that its size should be divided tenfold, it will nevertheless be 
%% %% published someday in some more reasonable journal.
%% %% \footnote{the author hadalready been subjected to similar 
%% %% inspired comments on works such as 
%% %% \cite{conductor} and as a result this
%% %% and a number of other pieces never were 
%% %% published, and others have been
%% %% but after many years of delay.}  
%% %% 
%

The formula obtained in \cite{cras2002} is:
\begin{theorem}
Let $\phi_a^\pm$ be the even entire functions solving:
\begin{subequations}
\begin{align}\label{eq:phiplus}
    \phi_a^+(x) + \int_0^a 2\cos(2\pi xy)\phi_a^+(y)\,dy &= 2\cos(2\pi ax)\\
\label{eq:phimoins}
    \phi_a^-(x) - \int_0^a 2\cos(2\pi xy)\phi_a^-(y)\,dy &= 2\cos(2\pi ax)
\end{align}
\end{subequations}  
The function $\wh{E_a}(s)$ given as 
\begin{equation}
\wh{E_a}(s) = \sqrt{a}\left(a^{-s}  + \frac12\int_a^\infty (\phi_a^+(x) -
      \phi_a^-(x))x^{-s}\,dx\right)
\end{equation}
is an entire function and is such that for all $z,w\in\CC$
\begin{equation}
  \label{eq:repro2}
 \wh{X_z^a}(w) = \wh{X_w^a}(z) =  \int_a^\infty X_z^a(x)X_w^a(x)\,dx =
    \frac{\wh{E_a}(z) \wh{E_a}(w) - 
  \wh{\cC(E_a)}(z) \wh{\cC(E_a)}(w)}{z + w -1}
\end{equation}  
\end{theorem}

Equation \eqref{eq:repro2}  looks different from \eqref{eq:repro1} because we
dropped Gamma factors and also because we are not using the Hilbert scalar
product to avoid anti-analyticity.  Equations \eqref{eq:phiplus},
\eqref{eq:phimoins} are first considered in $L^2(0,a;dx)$. They have solutions
because $\pm1$ are not in the spectrum of the finite cosine kernels. They
imply that $\phi_a^+(x)$ and $\phi_a^-(x)$ are even entire functions, in
particular they have meaning for $x>a$. The integral for $\wh{E_a}(s)$ is
absolutely convergent for $\Re(s)>0$.

The function $\cE_a(s)$ is real for $s$ real. Writing $\cE_a(s) = \cA_a(s) -
i\cB_a(s)$ on the critical line, with $\cA_a(s)$ and $\cB_a(s)$ real there, we
thus have for $s$ on the critical line: $2\cA_a(s) = \cE_a(s)+
\overline{\cE_a(s)} = \cE_a(s) + \cE_a(\overline s) = \cE_a(s) +
\cE_a(1-s)$. The last equation, for general $s$, defines $\cA_a$ as an entire
function, which is even under $s\mapsto 1-s$, real for $s$ real, and real for
$s$ on the critical line. Also $2\cB_a(s) = i(\cE_a(s)-\cE_a(1-s))$ is an
entire function, which is odd under  $s\mapsto 1-s$, imaginary for $s$ real,
and real for $s$ on the critical line. We shall write $\cA_a(s) = \pi^{-\frac
s2}\Gamma(\frac s2)\wh{A_a}(s)$, $\cB_a(s) = \pi^{-\frac s2}\Gamma(\frac
s2)\wh{B_a}(s)$, so that $A_a(x)$ (resp. $B_a(x)$) is a cosine invariant
(resp. skew-invariant) distribution supported in $[a,+\infty)$. With
$\delta_a(x) = \delta(x-a)$ one has indeed:
\begin{subequations}
\begin{align}
    A_a(x) &= \frac{\sqrt a}2\left(\delta_a(x) + \phi_a^+(x)\Un_{x>a}(x)\right)\\
  -iB_a(x) &= \frac{\sqrt a}2\left(\delta_a(x) - \phi_a^-(x)\Un_{x>a}(x)\right)
\end{align}
\end{subequations} 
The evaluator identity \eqref{eq:repro1} can also be written as 
\begin{equation}
  \label{eq:repro4}
    (\cZ_w|\cZ_z) = 2\frac{\overline{\cA_a(z)}\cB_a(w) -\overline{\cB_a(z)}\cA_a(w)}{i(\overline z + w -1)} 
\end{equation}
In particular, for $z=w$ one obtains
\begin{align}
\label{eq:evalnorm}
  \|\cZ_z\|^2 &= 2\frac{\Im(\overline{\cA_a(z)}\cB_a(z))}{\Re(z)-\frac12}\qquad(\Re(z)\neq\frac12)
\\ \label{eq:evalnormcrit}
  \|\cZ_{\frac12+iE}\|^2 &= -2{\cA_a(\frac12+iE)}\,i\cB_a'(\frac12+iE) + 2{\cB_a(\frac12+iE)}\,i\cA_a'(\frac12+iE)
\end{align}     Derivatives in \eqref{eq:evalnormcrit} are taken with respect
to $s = \frac12+iE$ (not $E$). Equation \eqref{eq:evalnorm} implies that the
entire functions $\cA_a$ and $\cB_a$ have all their zeros on the critical
line. It is indeed a fact that all vectors $\cZ_z$ are non-zero (it is easy to
see that ``$\cF(z)=0$ for all $\cF$'' is impossible).  From
\eqref{eq:evalnorm} we deduce that the ratio $\frac{\cB_a(s)}{\cA_a(s)}$ is of
positive imaginary part when $\Re(s)>\frac12$. It is elementary, as
$\frac{\cB_a(s)}{\cA_a(s)}$  is meromorphic, that  this implies that its poles
and zeros (all on the critical line) are simple. Next, if  $\cA_a$ had a
multiple zero $\rho$, then $\cB_a$ would have to vanish there, but then taking
a limit in \eqref{eq:evalnorm} we would obtain   $\|\cZ_\rho\| = 0$, which
we have already mentioned is false. So the zeros of $\cA_a$ and $\cB_a$ are
simple; using \eqref{eq:evalnormcrit} one argues that the respective zeros
interlace. These arguments apply in the general situation of \cite{Bra}. A
simplification which is valid here and that we have made use of is that the
evaluators $\cZ_z$ are non-zero also when $z$ is on the symmetry line.

From \eqref{eq:repro4} two vectors $\cZ_{\rho}$ and $\cZ_{\rho'}$ associated
to distinct zeros of $\cA_a$ are always orthogonal. Furthermore from
\eqref{eq:repro4} if $\rho_n = \frac12+iE_n$ is a zero of $\cA_a$ then
$\cZ_{\rho_n}(s) = 2{i\cB_a(\rho_n)}\frac{\cA_a(s)}{s-\rho_n}$ (we used $\rho_n = 1
-\overline{\rho_n}$). So the vectors $\frac{\cA_a(s)}{s-\rho_n}$ are mutually
orthogonal. From \eqref{eq:evalnormcrit}:
\begin{equation}
 \left\|\frac{\cA_a(s)}{s-\rho_n}\right\|^2_{K_a} = \int_{\Re(s)=\frac12}
\left|\frac{\wh{A_a}(s)}{s-\rho_n}\right|^2\,\frac{|ds|}{2\pi} = \frac12
\left.\frac{d}{dE}\right|_{E=E_n} \frac{\cA_a(\frac12+iE)}{\cB_a(\frac12+iE)}
\end{equation}
We can write this also as
\begin{equation}
  \left\|\frac{\cA_a(s)}{s-\rho_n}\right\|^2_{K_a} =   = \frac12
\left.\frac{d}{dE}\right|_{E=E_n}
\frac{\wh{A_a}(\frac12+iE)}{\wh{B_a}(\frac12+iE)} = \frac1{2} \left.\frac{d}{i\,dE}\right|_{E=E_n}
\frac{\wh{A_a}(\frac12+iE)}{\wh{E_a}(\frac12+iE)}
\end{equation} The vectors $\frac{\cB_a(s)}{s-\rho}$  corresponding to the
zeros of $\cB_a$ provide another orthogonal system.  
These vectors corresponding to the zeros of $\cA_a$ (or $\cB_a$) span all of
$K_a$. This orthogonal basis theorem is a general fact from \cite{Bra} (under
circumstances which are verified here).  

In the next sections we prove this orthogonal basis theorem from the standard
material provided  by the Weyl-Stone-Titchmarsh-Kodaira theory of second order
differential equations (\cite{coddlevinson,levsarg}), as an outcome of
studying the deformation with respect to $a>0$ of the spaces $K_a$
(\cite{cras2003}). Anyhow, obviously one of the main origins of the axioms of
\cite{Bra} is precisely the behavior of the objects from the Weyl-Stieltjes
spectral theory of differential or difference equations  (see Remling
\cite{remling} for the spectral problems of Schrödinger equations from the
point of view of these axioms).

\section{Dirac  system and associated isometric expansion}

\begin{definition}
  We let $C_a$ be the kernel $2\cos(2\pi xy)$ acting on
  $L^2(0,a;dx)$. And we define
\begin{equation}
 \mu(u) = a\frac{d}{da}\log\frac{\det(1+C_a)}{\det(1-C_a)}\qquad\text{with\
} u = \log(a)
\end{equation}
\end{definition}

\begin{definition}
  We let $\cA(u,s) = \cA_a(s)$ and $\cB(u,s) = \cB_a(s)$ for $a = \exp(u)$,
  $-\infty<u<\infty$. 
\end{definition}

\begin{theorem}[\cite{cras2003}]
  The following differential system holds, where $s= \frac12+iE$, $E\in\CC$:
  \begin{subequations}
  \begin{align}
    (\frac{d}{du}   + \mu(u)) \cA(u,s) &= - E\cB(u,s)\\
    (\frac{d}{du}   - \mu(u)) \cB(u,s) &= + E\cA(u,s)
  \end{align}
  \end{subequations}
Equivalently:
  \begin{equation}
    \label{eq:36}
    \left(\begin{bmatrix}
      0&1\\-1&0\end{bmatrix}\frac{d}{du} - 
    \begin{bmatrix} 0&\mu(u)\\\mu(u)&0\end{bmatrix}\right)
      \begin{bmatrix} \cA(u,s)\\  \cB(u,s) \end{bmatrix} 
= E \begin{bmatrix} \cA(u,s)\\ \cB(u,s) \end{bmatrix}
  \end{equation}
\end{theorem}

Combining the differential system \eqref{eq:36} and the evaluator formula
\eqref{eq:repro4} one obtains the following identity
\begin{equation}\label{eq:repro5} (\cZ_w^a|\cZ_z^a) = \int_{\log(a)}^\infty
\overline{2\cA(u,z)}\,2\cA(u,w) + \overline{2\cB(u,z)}\,2\cB(u,w)\; \frac{du}2
\end{equation}  and in particular we see that  for each $z\in\CC$ the
two-component vector $\begin{bmatrix} \cA(u,z)\\  \cB(u,z) \end{bmatrix}$  is
square-integrable at $+\infty$.

Let us recall that according to a general theorem \cite[\textsection13, Thm
7.1]{levsarg} for Dirac systems with continuous coefficients we are
automatically in the limit-point case at infinity. Being in the limit point
case at $+\infty$  means that for each $E\not\in\RR$ there is exactly one (up
to a scalar multiple) non-vanishing solution which is square-integrable at
$+\infty$. Thus we have identified this solution, and we have seen that it
also exists for $E$ real (the solutions which are square-integrable at
$-\infty$ are also known \cite{cras2003}; we will mention them in a later
section). There can't be other linearly independent square-integrable
solutions as this would characterize the limit-circle case, which never
happens for Dirac type systems with continuous coefficients.\footnote{in
\cite{hankel} we give an eigensolution for $E=0$ which we prove to be not
square-integrable at infinity.  This explicit argument confirms that we are
not in the limit-circle case.}

Let us consider the Dirac-type differential operator
$\displaystyle\begin{bmatrix} 0&1\\-1&0\end{bmatrix}\frac{d}{du} - 
    \begin{bmatrix} 0&\mu(u)\\\mu(u)&0\end{bmatrix}$ acting on pairs
$\begin{bmatrix} \alpha(u)\\  \beta(u) \end{bmatrix}$, $-\infty<u<+\infty$. As
we are in the limit-point case at infinity (plus or minus), from Hermann
Weyl's theory (adapted to first-order two-component systems) we know that the
symmetric operator defined by this differential operator with domain given by
the functions of class $C^1$ and with compact support is essentially
self-adjoint (\cite{levsarg,coddlevinson,reedsimon}). In Theorem 15 of \cite{hankel} we determined the associated
isometric spectral expansion and showed that this self-adjoint operator $H_0$
is $-i(x\frac{d}{dx}+\frac12)$ acting on $L^2(0,\infty;dx)$
or equivalently the operator of multiplication by $E$ on the critical
line:

\begin{theorem}[\cite{hankel}, Theorem 15]\label{thm:expansion}
There are  unitary identifications between $L^2(0,\infty;dx)$,
$L^2(\Re(s)=\frac12;\frac{|ds|}{2\pi|\pi^{-\frac
        s2}\Gamma(\frac s2)|^2})$ and $L^2(\RR\to \CC^2;\frac{du}2)$
given  by the following formulas (to be understood as suitable $L^2$-limits,
see \cite{hankel}), where always $s= \frac12+iE$, $E\in\RR$:
\begin{align}
\cF(s) &= \pi^{-\frac s2}\Gamma(\frac s2)\int_{0<x<\infty} f(x){x^{-s}}\,dx \\
\label{eq:iso1}
\begin{bmatrix} \alpha(u)\\  \beta(u) \end{bmatrix} &= \int_{E\in\RR}
    \cF(s)\,
    {\begin{bmatrix}
      {2\cA(u,s)}\\ {2\cB(u,s)}
    \end{bmatrix}}
  \;\frac{dE}{2\pi\, |\pi^{-\frac
        s2}\Gamma(\frac s2)|^2}\\
\label{eq:iso2}
\cF(s) &=  \int_{u\in\RR} \alpha(u)\,2\cA(u,s) +
  \beta(u)\,2\cB(u,s)\,\frac{du}2 \\
\int_0^\infty |f(x)|^2\,dx &= \int_{-\infty}^\infty |\cF(s)|^2 \;\frac{dE}{2\pi\, |\pi^{-\frac
        s2}\Gamma(\frac s2)|^2} = \int_{-\infty}^\infty (|\alpha(u)|^2 +
|\beta(u)|^2)\frac{du}2
\end{align}
Under these identifications:
\begin{itemize}
\item The orthogonal projection of $f(x)$ to $K_{a_0}$
corresponds to $
\left[\begin{smallmatrix}
  \alpha(u)\\ \beta(u)
\end{smallmatrix}\right] \mapsto \left[\begin{smallmatrix}
  \alpha(u)\\ \beta(u)
\end{smallmatrix}\right] \Un_{u>\log(a_0)}(u)$,
\item  the Cosine Transform $f\mapsto \cC(f)$ corresponds  to
$
\left[\begin{smallmatrix}
  \alpha\\ \beta
\end{smallmatrix}\right] \mapsto \left[\begin{smallmatrix}
  \alpha\\ -\beta
\end{smallmatrix}\right]$ (and $\cF(s)\mapsto\cF(1-s)$),
\item the
self-adjoint operator $-i(x\frac{d}{dx} + \frac12)$ on $L^2(0,+\infty;dx)$,
which multiplies by $E$ with maximal domain on the critical line,
corresponds to the unique self-adjoint completion of the differential operator on the $u$-line: 
  \begin{equation}
    H_0 = \begin{bmatrix}
      0&\frac{d}{du}\\-\frac{d}{du}&0\end{bmatrix} - 
    \begin{bmatrix} 0&\mu(u)\\\mu(u)&0\end{bmatrix}
  \end{equation} 
with initial domain 
the functions of class $C^1$ (or even $C^\infty$) with
compact support.
\end{itemize}
\end{theorem}

This Theorem is an elaboration deduced mainly from the evaluator identity
\eqref{eq:repro5}. Let us also mention that the integrals in \eqref{eq:iso1}
are absolutely convergent when $f$ is in the domain of $H_0$, that is, when
$s\cdot \wh{f}(s)$ is square integrable on the critical line. The functions
$\alpha(u)$ and $\beta(u)$ are then continuous; actually we know that they are
absolutely continuous as they must be in the domain of $H_0$ in the
$u$-picture, so their derivatives must be locally square integrable.

The operator $xp = -ix \frac{d}{dx}$, or $\frac12(xp+px) =  -i(x
\frac{d}{dx}+\frac12)$ was considered ``semi-classically'' by Berry
and Keating \cite{BK1} in a  discussion where it was accompanied with
the conditions $|x|>l_x$, $|p|>l_p$ (the operator $-i(x
\frac{d}{dx}+\frac12)$ with conditions $|x|,|p|<\Lambda$ arose in the
work of Connes \cite{Connes}). Certainly it is natural to consider
that the conditions $|x|>l_x$, $|p|>l_p$,  treated
``semi-classically'' in \cite{BK1}, could at the ``quantum'' (that is,
operator-theoretical level) be related with the condition of support
which defines the spaces $K_a$  (and if we were considering functions
odd in $x$, we would have the analogues of the spaces $K_a$ associated
with the sine kernel). Recently, German Sierra in \cite{Sierra4} (see
also \cite{Sierra2,Sierra3}) has generalized  the conditions
$|x|>l_x$, $|p|>l_p$ through the consideration of
differential-integral equations obtained from the differential
operator $xp$   by the addition of (singular) finite rank
perturbations. Sierra presents general aspects of the corresponding
eigenfunctions. This leads to scattering states (deformations of the
``free'' states $x^{iE} x^{-\frac12}$), and possibly also to bound
states in some fine-tuned situations. In the next section we see how
we have natural bound states in $K_a$ and then we will show that in
the $x$-picture these states are eigensolutions in the form considered
in \cite[Section IV]{Sierra4}. The scattering states will span the
perpendicular complement of $K_a$ in $L^2$.

\section{Bound states}

We consider the differential system \eqref{eq:36} on the half-interval
$(u_0,+\infty)$, $u_0 = \log(a_0)$ with boundary condition $\alpha(u_0) =
0$. From general theorems \cite{reedsimon,levsarg,coddlevinson} we know that
this leads to a self-adjoint operator and associated isometric expansion (with
discrete spectrum from the behavior of $\mu(u)$ at $+\infty$, but we will confirm this directly). To obtain the spectrum we apply Hermann Weyl's method
(\cite[\textsection9]{coddlevinson}, \cite[\textsection3]{levsarg}). For each
$s\in\CC$ we let $\psi(u,s)$ be the unique solution of the system
\eqref{eq:36} for the eigenvalue $E$, $s=\frac12+iE$, and with the initial
condition $\psi(u_0,s) = \left[\begin{smallmatrix} 1\\ 0
\end{smallmatrix}\right]$ and let $\phi(u,s)$ be the unique solution for the
eigenvalue $E$ with the initial condition $\phi(u_0,s) =
\left[\begin{smallmatrix} 0\\ 1
\end{smallmatrix}\right]$. Let, for
$\Im(E)>0$, the quantity $m(E)$ be chosen such that
$\psi(u,s)-m(E)\phi(u,s)$ is the 
unique linear combination which has the property of being square-integrable at
$+\infty$. 

To check up on signs, let us see what would happen if the potential $\mu(u)$
were identically vanishing. Then we  would have $\psi(u,s) =
\left[\begin{smallmatrix} \cos(E(u-u_0))\\
\sin(E(u-u_0)) \end{smallmatrix}\right]$ and $\phi(u,s) =
\left[\begin{smallmatrix} -\sin(E(u-u_0))\\
\cos(E(u-u_0)) \end{smallmatrix}\right]$, and the $m$ function would be the
constant $+i$. In the general case the $m$ function is known to be an analytic
function mapping the upper half-plane to itself. Here, as the one-dimensional
space of square-integrable solutions has basis $\left[\begin{smallmatrix}
\cA(u,s)\\ \cB(u,s)
\end{smallmatrix}\right]$  we have $m(E) =
\frac{-\cB_{a_0}(s)}{\cA_{a_0}(s)}$. We had already argued from the evaluator
formula \eqref{eq:repro4} that this meromorphic function has positive
imaginary part on the half-plane $\Im(E)>0$. As our  $m$-function is
meromorphic and real-valued for $E$ real, the spectral measure $\nu$, which is
given  via the formula $\nu(a,b) = \lim_{\epsilon\to 0^+} \frac1\pi \int_a^b
\Im( m(E+i\epsilon))\,dE$ (under the condition $\nu\{a\} +\nu\{b\} = 0$) is  a
sum of Dirac measures and the spectrum is thus purely discrete, having its
support on the zeros $E_n$ of $\cA_{a_0}({\frac12+iE})$. We
obtain:\footnote{derivatives are computed with respect to $s = \frac12+iE$ so
$i\cA_{a_0}'(\frac12+iE_n)$ is real.}
\begin{equation}
  \label{eq:191}
  \nu(dE) = \sum_n
\frac{\cB_{a_0}(\frac12+iE_n)}{i\cA_{a_0}'(\frac12+iE_n)}\delta(E-E_n)\;dE
\end{equation}
From the general theory the Hilbert space
$L^2((u_0,+\infty)\to \CC^2;du)$ admits an orthogonal basis $(T_n)$
given by the two-component vectors
$T_n(u)
= \phi(u,\frac12+iE_n)$
which have, according to the spectral measure formula \eqref{eq:191}, norms 
$\int_{u_0}^\infty |T_n(u)|^2\,du  =
\frac{i\cA_{a_0}'(\frac12+iE_n)}{\cB_{a_0}(\frac12+iE_n)}$. The
spectral measure is $ \nu(dE) = \sum_n \frac{\delta(E-E_n)}{\|T_n\|^2} \,dE$ and
the Plancherel formula is
\begin{equation}
\int_{u_0}^\infty |T(u)|^2\,du  = \sum_n \frac{|\wt{T}(E_n)|^2}{\int_{u_0}^\infty |T_n(u)|^2\,du }\qquad \wt{T}(E_n) =
\int_{u_0}^\infty T(u)\cdot\phi(u,\frac12+iE_n)\,du
\end{equation}
Let us consider the vectors
\begin{equation}
 Z_n(u) = \left[\begin{smallmatrix} 2\cA(u,\frac12+iE_n)\\
2\cB(u,\frac12+iE_n)\;
\end{smallmatrix}\right]\Un_{u>u_0}(u) = 2\cB(u_0,\frac12+iE_n)\;T_n(u)
\end{equation} We deduce from what has just been
stated regarding the $T_n$'s that the $Z_n$'s are an orthogonal basis of
$L^2((u_0,+\infty)\to \CC^2;\frac12du)$ with Hilbert norms in this space:
\begin{equation}\label{eq:Zn}
 (Z_n|Z_n) =  2 i \cA_{a_0}'(\frac12+iE_n)\cB_{a_0}(\frac12+iE_n)
\end{equation} From Theorem \ref{thm:expansion} we know that $Z_n$ in
the ``$\cF(s)$-picture'' corresponds precisely to the evaluator
$\cZ_{\frac12+iE_n}^{a_0}$. Equations \eqref{eq:Zn} and
\eqref{eq:evalnormcrit} do indeed match. Furthermore the two-component
vector $Z_n(u)$ solves the differential system on $(u_0,+\infty)$ for
the eigenvalue $E_n$. In conclusion the differential operator
$\displaystyle\begin{bmatrix} 0&1\\-1&0\end{bmatrix}\frac{d}{du} -
    \begin{bmatrix} 0&\mu(u)\\\mu(u)&0\end{bmatrix}$ acting on pairs
$\left[\begin{smallmatrix} \alpha(u)\\ \beta(u) \end{smallmatrix}\right]$  on
$(u_0,\infty)$ with the boundary condition $\alpha(u_0) = 0$ is self-adjoint
with a purely discrete spectrum given by the imaginary parts $E_n$ of the
zeros of $\cA_{a_0}$ ($a_0 = \exp(u_0)$). We notice, as
$\cZ_{\frac12+iE_n}^{a_0}(\frac12+iE) = 2\cB_{a_0}(\frac12+iE_n)
\frac{\cA_{a_0}(\frac12+iE)}{E - E_n}$, that the eigenvectors
$T_n(u)$ which are such that $T_n(u_0) =
\left[\begin{smallmatrix}
  0\\1
\end{smallmatrix}\right]$ correspond exactly in the $\cF(s)$-picture to the
functions $s\mapsto \frac{\cA_{a_0}(\frac12+iE)}{E - E_n}$.

Let us denote by $H_0^{a_0}$ the self-adjoint operator given by the differential
system \eqref{eq:36} on $(u_0,+\infty)$ with boundary condition
$\left[
  \begin{smallmatrix}
    0\\ *
  \end{smallmatrix}\right]$ at $u_0 = \log(a_0)$. With reference to the comment after
Theorem \ref{thm:expansion}, we observe that its eigenvectors (in the $\cF(s)$
picture) $s\mapsto \frac{\cA_{a_0}(s)}{s -\rho}$ ($\cA_{a_0}(\rho) = 0$) are
not located in the domain of $H_0$ (which multiplies by $E$), and indeed
$\wh{A_{a_0}}(s)$ is not square-integrable on the critical line, as can be
seen directly or from the presence of the Dirac at $a$ in $A_{_0}(x)$ or the
behavior of $A_{a_0}(x)$ as $x\to +\infty$. To discuss this further let us
mention:

\begin{proposition} If $f$ is in $K_a$ and is in the domain of $H_0$ in $L^2$,
then $H_0(f)$ is in $K_a$. The restriction of $H_0$ to $K_a$ has deficiency
indices $(1,1)$. The self-adjoint operator $H_0^a$ is an extension of this
restriction.
\end{proposition}

Let  $f$ in $K_a$ be such that $H_0(f)$ is in $L^2$. In the $u$-picture the
two-component vector $T$ representing $f$ must have its components absolutely
continuous, and the action of $H_0$ on $T$ is the one given as a differential
operator. As $f$ is in $K_a$, the components of $T$ vanish identically for
$u<\log(a)$ so their derivatives also. Hence $H_0(f)\in K_a$. The statement is
also easy to establish directly  in the $x$-picture using the formula $\cC
(x\frac{d}{dx} + \frac12) f  = - (x\frac{d}{dx} + \frac12) \cC(f)$ (in the
distribution sense). Let us note though that by necessity the function $T(u)$
being continuous will be such that $\alpha(\log(a))= 0$ and $\beta(\log(a)) =
0$. The domain of the restriction of $H_0$ to $K_a$, in the $u$-picture is
thus given by the absolutely continuous square integrable functions, with both
components vanishing
at $\log(a)$, and whose images under the differential operator are
square-integrable. This shows that the operator $H_0^a$ is indeed an extension
of this restriction. From well-known theorems \cite{reedsimon,levsarg}, the
deficiency indices of such a differential operator on a semi-infinite
interval, in the limit-point case at infinity, are $(1,1)$ with the
one-parameter family of possible boundary conditions at $\log(a)$ giving the
self-adjoint extensions. We can confirm directly in the $\cF(s)$-picture the
computation of the deficiency indices. Indeed $i(H_0 - \frac i2)$ is
multiplication by $s$. As $s\cF(s)$ vanishes at $0$ (here $\cF$ is taken in a
$K_a$, so is also defined away from critical line) it is perpendicular
(if $\cF$ is in the domain of $H_0$) to the evaluator $\cZ_0^a$. Conversely if
$\cG$ in $K_a$ is perpendicular to $\cZ_0^a$, the quotient $\frac1s \cG(s)$ is
an entire function and it is easy to prove that it corresponds to an element
of $K_a$ (see \cite{twosys} for details). So the range of
$\left.H_0\right|_{K_a} - \frac i2$ has codimension $1$, and similarly for the
range of $\left.H_0\right|_{K_a} + \frac i2$. This confirms that the deficiency
indices are $(1,1)$.

\section{Continuous spectrum}

\begin{definition} We let $L_a$ be the perpendicular complement of
$K_a$ in $L^2$. In the $x$ picture one has $L_a = L^2(0,a;dx) +
\cC(L^2(0,a;dx))$. In the $u$ picture $L_a$ is the space of the square
integrable two-component vectors supported on $(-\infty,\log(a))$.
\end{definition}

Again, from the $u$-picture we see that if $f$ is in $L_a$ and is in the
domain of $H_0$ in $L^2$ then
$H_0(f)$ is again in $L_a$. And we again have  a one-parameter family of
self-adjoint extensions of the restriction of $H_0$ to $L_a$. But all these
extensions will have continuous spectra, as the potential function $\mu(u)$
vanishes exponentially as $u\to-\infty$. Summing up:

\begin{proposition} If $f$ is in $L_a$ and is in the domain of $H_0$ in $L^2$,
then $H_0(f)$ is in $L_a$. The restriction of $H_0$ to $L_a$ has deficiency
indices $(1,1)$. Its self-adjoint extensions correspond to a one-parameter
family of boundary conditions at $\log(a)$ and each has a purely continuous
spectrum.
\end{proposition}

To obtain a self-adjoint operator on all of $L^2$ we thus need to choose the
boundary condition at $u_0$. A specific choice will emerge from the following
procedure. The explicit spectral expansion will be given in the next section,
after having explained how the results obtained in the present section fit
into the framework considered by Sierra in \cite{Sierra4}, once we go back to
the $x$-picture.

We consider $H_0$ as a differential operator in the $u$ picture. The bound
states $\frac{\cA_{a}(\frac12+iE)}{E - E_n}$ are there the
vectors denoted $T_n(u)$, which have support in $u\geq u_0$, and verify
$\lim_{u\to u_0} T_n(u) = \left[
  \begin{smallmatrix}
    0\\1
  \end{smallmatrix}\right]$. They obey for $-\infty<u<+\infty$ the
differential equation, say as distributions:
\begin{equation}
  \begin{bmatrix}
      0&\frac{d}{du}\\-\frac{d}{du}&0\end{bmatrix}  T_n(u)- 
    \begin{bmatrix} 0&\mu(u)\\\mu(u)&0\end{bmatrix}  T_n(u) = E_n T_n(u) +\begin{bmatrix} \delta(u-u_0)\\ 0\end{bmatrix}
\end{equation}
We will write (recalling the use of $\frac{du}2$ in the $u$ picture):
\begin{equation}
 \delta_0(u) = \begin{bmatrix} 2\delta(u-u_0)\\ 0\end{bmatrix}\qquad
\delta_1(u) = \begin{bmatrix} 0\\2\delta(u-u_0)\end{bmatrix}
\end{equation} 
and shall convene that $\left[
\begin{smallmatrix}
  (\delta_0 | T)\\(\delta_1 | T)
\end{smallmatrix}\right] = \frac12(\lim_{u\to u_0^+} + \lim_{u\to u_0^-})T(u)$
when the limits exist. So here   $(\delta_0 | T_n) = 0$ and $(\delta_1 |T_n) =
\frac12$. This transforms the inhomogeneous equation into a linear homogeneous
equation:
\begin{equation}
  \begin{bmatrix}
      0&\frac{d}{du}\\-\frac{d}{du}&0\end{bmatrix}  T_n(u)- 
    \begin{bmatrix} 0&\mu(u)\\\mu(u)&0\end{bmatrix}  T_n(u) = E_n T_n(u) +
(\delta_1 |T_n) \delta_0(u)
\end{equation} In order to obtain a formally symmetric expression
we add one more term, obtaining in the end the following perturbation of the differential
operator $H_0$:
\begin{equation}\label{eq:Hoperator}
 H(T) = H_0(T) - (\delta_1 | T) \delta_0 - (\delta_0 | T) \delta_1 
\end{equation}
According to the remark made after Theorem \ref{thm:expansion}, as a corollary
to the fact that  the functions $\frac{\wh{A_a}(s)}s$ and
$\frac{\wh{B_a}(s)}s$ are square-integrable on the critical line, when $T$ is
in the domain of $H_0$ its components are continuous functions and evaluation
at $u_0$ is a continuous linear form (given by equation \eqref{eq:iso1}). As
a bilinear form, the perturbation term is thus bounded relative to $(T|T) +
(H_0(T)|H_0(T))$.

Let us consider for a given $E$ the equation $H(T) = ET$ on the full
line. So $T$ must be of the form $T(u) = L(u)\Un_{u<u_0}(u) +
R(u)\Un_{u>u_0}(u)$ where $L$ and $R$ verify on the full line (as
differential equations) $H_0(L) = E L$, $H_0(R) = E R$. Writing
$L(u_0) = \left[
\begin{smallmatrix}
  l_0\\l_1
\end{smallmatrix}\right]$ and  $R(u_0) = \left[
\begin{smallmatrix}
  r_0\\r_1
\end{smallmatrix}\right]$ we obtain from \eqref{eq:Hoperator} the
matching conditions:
\begin{equation}
 -l_1 \delta_0 + l_0 \delta_1 + r_1 \delta_0 - r_0 \delta_1 = (l_1+r_1)
\delta_0 + (l_0+r_0) \delta_1 \quad \iff l_1 = r_0 = 0 
\end{equation}
We thus have for each $E$ a two dimensional space of eigensolutions,
which are of the form:
\begin{equation}
 T_E(u) = \lambda L_E(u)\Un_{u<u_0}(u) + \mu R_E(u)\Un_{u>u_0}(u)
\end{equation}
where $L_E$ verifies the differential equation $H_0(L_E) = E L_E$ and the
condition $L_E(u_0) = \left[
\begin{smallmatrix}
 1\\0
\end{smallmatrix}\right]$ and $R_E$ verifies the differential equation $H_0(R_E) = E R_E$ and the
condition $R_E(u_0) = \left[
\begin{smallmatrix}
 0\\1
\end{smallmatrix}\right]$.

As we saw in the last section the vectors $R_E(u)\Un_{u>u_0}(u)$ lead to a
well-defined spectral expansion of $K_a$ where in fact only those among the
$R_E$ which are square-integrable contribute, the suitable $E$'s being the
(imaginary parts of the) zeros of the entire function $\cA_{a_0}(s)$. On the
other hand, on the basis of standard theorems of
spectral expansions \cite{reedsimon, levsarg, coddlevinson} the vectors
$L_E(u)\Un_{u<u_0}(u)$ for $E$ real are the generalized eigenvectors of the
self-adjoint extension of the restriction of $H_0$ to $(-\infty,\log(a))$,
given by the boundary condition $\left[
\begin{smallmatrix}
 * \\ 0
\end{smallmatrix}\right]$ at $u_0$.

We now extend $H_0^{a}$ to all
of $L^2$ defining it on $L_{a}$ to be this self-adjoint extension of
$\left.H_0\right|_{L_{a}}$. In the  end we have thus been led by this
procedure to a well-defined 
self-adjoint operator on $L^2$ whose spectrum has on one hand a discrete part
corresponding to $K_a$ and on the other hand a continuous part corresponding
to $L_a$. Its eigenvectors are the solutions of the equations:
\begin{equation}
 H_0(T) - (\delta_1 | T) \delta_0 - (\delta_0 | T) \delta_1  = E T
\end{equation}
But only those solutions that are square-integrable as $u\to+\infty$ have to be
retained, those solutions being: on one hand the vectors previously denoted
$T_n(u)$ on $(u_0,+\infty)$, on the other hand a continuous spectrum of
generalized elements living on $(-\infty,u_0)$.
The spectral density will be determined in the next section. Here we now
explain how the $x$ picture looks like.

According to Theorem \ref{thm:expansion}, equation \eqref{eq:iso1}, 
$\delta_0(u)$ corresponds in the
$\cF(s)$ picture to ${2\cA_{a_0}(s)}$ and $\delta_1(u)$ corresponds to
${2\cB_{a_0}(s)}$ (here we are strictly on the critical line, these functions
are real-valued).  In the
$\wh{f}(s)$ picture,  $\delta_0(u)$ corresponds thus to $2\wh{A_{a}}(s)$ and
$\delta_1(u)$ to $2\wh{B_{a}}(s)$. In the $x$ picture $\delta_0$ and
$\delta_1$ correspond to the distribution $2A_a(x)$ and $2B_a(x)$. The
operator $H$ takes the form:
\begin{equation}
  \label{eq:Hx}
  H(f) = -i(x\frac{d}{dx} + \frac12)f - (2B_a|f) 2A_a - (2A_a|f)2B_a\,,
\end{equation} where the ``scalar products'' $(2A_a|f)$ and $(2B_a|f)$ have to
be defined in a suitable manner if $f$ does not belong to the domain of
$H_0$. We will switch to another linear combination allowing an easier
discussion and having the format considered by Sierra in \cite{Sierra4}.

Using $E_a = A_a - i B_a$ and $F_a = \cC(E_a) = A_a + i B_a$ we
rewrite equation \eqref{eq:Hx} as 
\begin{equation}
  \label{eq:Hx2}
  H(f) = -i(x\frac{d}{dx} + \frac12)f  + 2 i(E_a|f)F_a - 2i(F_a|f)E_a
\end{equation} The (real valued) distribution $E_a$ differs from a
Dirac $\sqrt{a}\,\delta_a(x)$ by a square integrable function
supported on $[a,+\infty)$, and its Cosine transform $F_a = \cC(E_a) =
A_a + i B_a$ is a function supported on $[a,+\infty)$  which differs
from $\sqrt{a}\, 2\cos(2\pi ax)$ by $O(\frac1x)$ at $+\infty$.  When
we try, directly in the $x$ picture, to solve the eigen-equation $H(f)
= Ef$ we have to give a meaning to $(E_a|f)$ where $f$ may have a jump
at $a$ (we may give an imaginary part to $E$ and look for a
square-integrable $f$ to simplify the discussion). The correct choice
is to take $(\delta_a|f)$ for those $f$ to be the average between the
left and right values at $x=a$. And $(F_a|f)$ needs to be defined
suitably as well. The discussion in the $u$-picture was technically
easier and we shall defer to another publication the more detailed
direct treatment of \eqref{eq:Hx2} in the $x$-picture.

Let us mention that if we repeat the whole procedure but starting from the
self-adjoint operator on $K_a$ corresponding to the zeros of the $\cB_a$
function we end up with the following operator
\begin{equation}
  \label{eq:Hx3}
  H(f) = -i(x\frac{d}{dx} + \frac12)f + 2i(F_a|f)E_a - 2 i(E_a|f)F_a 
\end{equation}
which differs from the previous one only by the sign of the perturbation. 

The general equation considered by Sierra in \cite[IV]{Sierra4} is
\begin{equation}
  \label{eq:Hx4}
  H(f) = -i(x\frac{d}{dx} + \frac12)f + i(\psi_2|f)\psi_1 - i(\psi_1|f)\psi_2
\end{equation} Thus \eqref{eq:Hx2} and \eqref{eq:Hx3} are indeed special
cases. Furthermore Sierra considers  the ``duality'' condition that $\psi_1$
and $\psi_2$ are a Cosine (or Sine) transform pair. Our special instances do
have this property. Finally in \cite[VI]{Sierra4}  the case $\psi_2(x) =
\sqrt2 \sqrt{a}\; \delta_a(x)$ and $\psi_1 = \cC(\psi_2)$ is devoted special
attention. Whether $\psi_2(x) = \sqrt2 \sqrt{a}\; \delta_a(x)$ leads  to bound
states is not clear.

\section{More on the scattering spectrum}

We want to provide more information on the ``scattering'' part of the
spectrum.  When $E$ is not real the differential system \eqref{eq:36} has (up
to a multiple) exactly one solution which is square-integrable at $-
\infty$. Fortunately, this solution is known explicitely. It was determined in
\cite{cras2003}, as well as the way $\left[\begin{smallmatrix} \cA(u,s)\\
\cB(u,s)
\end{smallmatrix}\right]$ for $E$ real is expressed as a linear combination of the
boundary values of such solutions and their complex conjugates.

\begin{theorem}[\cite{cras2003}; Theorem 16 of \cite{hankel}]
  The eigensolution $\begin{bmatrix} \alpha_E(u)\\  \beta_E(u) \end{bmatrix}$ of the differential system which is square-integrable at
  $-\infty$ for $\Im(E)>0$ ($\Re(s)<\frac12$) is provided by the following
  formulas (where $a = \exp(u)$):
  \begin{align}
   \alpha_E(u) = J(u,\frac12+iE) &= e^{-iEu} - e^{\frac12 u}\int_0^a \phi_a^+(x)
   x^{-s}\,dx\\ 
  \beta_E(u) = K(u,\frac12+iE) &= i e^{-iEu} + i e^{\frac12 u}\int_0^a
    \phi_a^-(x) x^{-s}\,dx
\end{align}
The following identity holds:
\begin{equation}\label{eq:scattering}
  \begin{bmatrix} 2\cA(u,\frac12+iE)\\  2\cB(u,\frac12+iE) \end{bmatrix} =
   \pi^{-\frac s2}\Gamma(\frac s2) \begin{bmatrix} J(u,\frac12+iE)\\
    K(u,\frac12+iE) \end{bmatrix} +  \pi^{-\frac{1-s}2}\Gamma(\frac{1-s}2)\begin{bmatrix} J(u,\frac12-iE)\\
    -K(u,\frac12-iE) \end{bmatrix}
\end{equation}
\end{theorem}

We will also use the notation $\wh{j_a}(s) = J(u,s)$ and
$\wh{k_a}(s) = K(u,s)$, for $u=\log(a)$:
\begin{subequations}
\begin{align}\label{eq:j}
  J(u,s) = \wh{j_a}(s) &= a^{\frac12-s} - \sqrt{a}\int_0^a \phi_a^+(x) x^{-s}\,dx\\
\label{eq:k}
  K(u,s) = \wh{k_a}(s) &= ia^{\frac12-s} + i \sqrt{a}\int_0^a \phi_a^-(x) x^{-s}\,dx
\end{align} 
\end{subequations} The functions $\phi_a^+(x)$ and $\phi_a^-(x)$ being entire
and even, we see that $\wh{j_a}$ and $\wh{k_a}$ are meromorphic in the entire
complex plane with possible poles at $s=1$, $s=3$, $s=5$, \dots Of course they
still provide after analytic continuation a solution of the differential
system with respect to $\log(a)$. We have
\begin{equation}
2 \cA_a(s) = \pi^{-\frac s2}\Gamma(\frac s2) \wh{j_a}(s)
+\pi^{-\frac{1-s}2}\Gamma(\frac{1-s}2) \wh{j_a}(1-s) 
\end{equation}  As $\cA_a(s)$ is
entire, if $\wh{j_a}$ does not have a pole at $1+2m$ ($m\in\NN$), necessarily the second
term on the right-hand side is forbidden to have a pole there, which means
that necessarily one must have a zero of $\wh{j_a}$ at $-2m$. We will now
argue that neither $\wh{j_a}$ nor $\wh{k_a}$ can have zeros in the half-plane
$\Re(s)<\frac12$ (which corresponds to $\Im(E)>0$). This will prove that
$\wh{j_a}$ (and similarly $\wh{k_a}$) does have poles at $s=1$, $s=3$, $s=5$,
\dots

Indeed as in the previous section we apply Hermann Weyl's theory
(\cite[\textsection9]{coddlevinson}, \cite[\textsection3]{levsarg}), but on
the interval $(-\infty,u_0)$. We maintain the same notations $\psi(u,s)$ with
boundary condition $\left[
  \begin{smallmatrix}
    1\\0
  \end{smallmatrix}\right]$ and 
$\phi(u,s)$  with
boundary condition $\left[
  \begin{smallmatrix}
    0\\1
  \end{smallmatrix}\right]$
but convene that the $m$-function for $\Im(E)>0$ is determined by
the condition that $-m(E)\psi(u,s)+\phi(u,s)$ is square integrable at
$-\infty$. Thus, we have the identity:
\begin{equation}
  \Im(E)>0\implies m(E) = \frac{-J(u_0,\frac12+iE)}{K(u_0,\frac12+iE)} = \frac{\wh{-j_{a_0}}(s)}{\wh{k_{a_0}}(s)}
\end{equation}
It is a general theorem that the $m$ function is an analytic function mapping
the upper half-plane to itself. So if $\wh{k_{a_0}}(s)$ has a zero for
$\Re(s)<\frac12$, then $\wh{j_{a_0}}(s)$ has to vanish also for the same $s_0$. But 
$\left[\begin{smallmatrix}
J(u,s_0)\\ K(u,s_0)
\end{smallmatrix}\right]$ is solution of a first order differential equation
with respect to $u$. If it vanishes for some $u$ it has to vanish for all
$u$'s. From \eqref{eq:j}, \eqref{eq:k} however it is easily argued that
$\wh{j_a}(s) \sim a^{\frac12-s}$ as $a\to 0^+$, and $\wh{k_a}(s)\sim_{a\to0}
ia^{\frac12-s}$ ($\Re(s)<1$).\footnote{the equations \eqref{eq:j}
  \eqref{eq:k} also give easily for a given $a$ asymptotic expansions of
  $\wh{j_{a}}(s)$, $\wh{k_a}(s)$ in powers of $\frac1s$ in the
  critical strip, and this gives thus corresponding information on $\wh{A_a}$,
  $\wh{B_a}$ and $\wh{E_a}$.}

The Wronskian $W(s) = \left|
  \begin{matrix}
    \cA(u,s)& J(u,s)\\
    \cB(u,s)& K(u,s)
  \end{matrix}\right|$ is independent of $u$. We can evaluate it at
$u\to-\infty$, perhaps first for $0<\Re(s)<\frac12$, using the asymptotic of
$J$ and $K$ as $a\to0^+$ and the identity \eqref{eq:scattering}. This
gives:
\begin{equation}
  \label{eq:wronskian}
  W(s) = i\,\pi^{-\frac{1-s}2}\Gamma(\frac{1-s}2)
\end{equation} Analytic continuation removes the restriction on $s$. This
confirms that $J$ or $K$ (in fact both as we saw) must have poles at $1$, $3$,
\dots

The Wronskian $W_1(u,s) = \frac{i}2\,\left|
  \begin{matrix}
    J(u,s)& \overline{J(u,s)}\\
    K(u,s)& \overline{K(u,s)}
  \end{matrix}\right| = \Im(-{J(u,s)}\overline{K(u,s)})$ verifies the differential equation
\begin{equation}\label{eq:diffwronskian}
 \frac{d}{du} W_1(u,s) = \Im(E)\, ( |J(u,s)|^2 + |K(u,s)|^2 )
\end{equation} For
$\Re(s)<1$, we have $W_1(u,s) \sim_{a\to0} a^{1-2\Re(s)}$. For $\Re(s) <
\frac12$, we have $\Im(E)>0$ and deduce from \eqref{eq:diffwronskian} that
$\Im(-{J(u,s)}\overline{K(u,s)})>0$ which confirms the $m$-property
of the ratio $\frac{-J(u,s)}{K(u,s)}$ (and the fact that none of $J$ and $K$
can vanish in this half-plane, whatever the value of $u$). Additionally, for
$\Re(s)=\frac12$ we obtain:
\begin{proposition}
  For $s= \frac12+iE$, $E\in\RR$ and $-\infty<u<+\infty$ there holds 
\begin{equation}
 \Im\left(-{J(u,s)}\,\overline{K(u,s)}\right) = 1
\end{equation}
Hence the functions $J(u,s)$ and $K(u,s)$ do not vanish on the line $\Re(s) = \frac12$.
\end{proposition}

Alternatively, we could also have deduced (not only for $E$ real but for all
$E\in\CC$)  the constant value $-2i$ of the Wronskian $ \left|\begin{matrix}
J(u,\frac12+iE)&  J(u,\frac12-iE)\\ K(u,\frac12+iE)& -K(u,\frac12-iE)
\end{matrix} \right|$ from the behavior at $-\infty$.  The
spectral measure $\nu(a,b) = \lim_{\epsilon\to 0^+} \frac1\pi \int_a^b \Im(
m(E+i\epsilon))\,dE$, $m(E) = \frac{-J(u_0,s)}{K(u_0,s)}$ is thus purely
absolutely continuous:
\begin{equation}
 \nu(dE) = \frac1{\pi\, |K(u_0,\frac12+iE)|^2} \;dE
\end{equation}
The spectral expansion takes the following form: we send $T(u) =
\left[\begin{smallmatrix}
  \alpha(u)\\ \beta(u)
\end{smallmatrix}\right]$ to 
\begin{equation}
  \label{eq:jTransform}
  \wt T(E) = \int_{-\infty}^{u_0} \alpha(u) {\psi_0(u,\frac12+iE)} +
\beta(u){\psi_1(u,\frac12+iE)} \,du
\end{equation}
where $\psi(u,s) = 
\left[\begin{smallmatrix}
  \psi_0(u,s)\\ \psi_1(u,s)
\end{smallmatrix}\right]$ is the unique solution of the differential system
verifying the boundary condition $
\left[\begin{smallmatrix}
  1\\ 0
\end{smallmatrix}\right]$ at $u_0$, \footnote{as $s=\frac12+iE$, $E\in\RR$, this
solution is real valued.} and this transform is isometric with the following
Plancherel identity:
\begin{equation}
 \int_{-\infty}^{u_0} |\alpha(u)|^2  + |\beta(u)|^2 \,\frac{du}2 = \frac1{2\pi}
\int_\RR |\wt T(E)|^2 \frac{dE}{|K(u_0,\frac12+iE)|^2} 
\end{equation}

The two linearly independent eigensolutions of the differential
system are $
\begin{bmatrix} J(u,\frac12+iE)\\ K(u,\frac12+iE) \end{bmatrix}$ and $
\begin{bmatrix} J(u,\frac12-iE)\\ -K(u,\frac12-iE) \end{bmatrix}$ with Wronskian $-2i$.  Hence:
\begin{equation}
  \psi(u,s) =  \frac{-i\,K(u_0,\frac12-iE) }{2}\begin{bmatrix} J(u,\frac12+iE)\\
    K(u,\frac12+iE) \end{bmatrix} + \frac{-i\,K(u_0,\frac12+iE) }{2}\begin{bmatrix} J(u,\frac12-iE)\\
    -K(u,\frac12-iE) \end{bmatrix}
\end{equation}
\begin{align}
 \psi_0(u, \frac12 + iE) &= \Im\left(K(u_0,\frac12-iE)J(u,\frac12+iE)\right)\\ 
\psi_1(u, \frac12 + iE) &= \Im\left(K(u_0,\frac12-iE)K(u,\frac12+iE)\right)
\end{align}

\end{document}